\documentclass[12pt]{amsart}
\usepackage{amsmath}
\usepackage{amssymb}
\usepackage{amsfonts}
\usepackage{tikz-cd}

\newtheorem{thm}{Theorem}

\newtheorem{lem}{Lemma}

\newtheorem{defn}{Definition}
\newtheorem{ex}{Example}

\begin{document}

\title{Inductive limits of quasi locally Baire spaces}
 \author{Thomas E. Gilsdorf}
 
\maketitle
\begin{center} Department of Mathematics \\
 Central Michigan University   \\
 Mt. Pleasant, MI  48859  \/ USA  \\ 
% ORCID:   0000-0001-9223-0388   \\
 \bigskip
 
  {\bf gilsd1te@cmich.edu}   \\
\bigskip
\today
 \end{center}
%   Email: \textbf{}.  

\bigskip 

\noindent  {\bf Abstract.}  \/ Quasi-locally complete locally convex spaces are generalized  to quasi-locally Baire  locally convex spaces.  It is shown that an inductive limit of strictly webbed spaces is regular if it is quasi-locally Baire.  This extends Qiu's theorem on  regularity.  Additionally, if each step is strictly webbed and quasi- locally Baire, then  the inductive limit is quasi-locally Baire if it is regular.   Distinguishing examples are provided.     \\
\bigskip



\bigskip

\section{Introduction and notation.}

\noindent  Inductive limits of locally convex spaces have  been studied in detail over many years. Such study includes properties that would imply regularity, that is, when every bounded subset in the in the inductive limit is contained in and bounded in one of the steps.  An excellent introduction to the theory of locally convex  inductive limits, including regularity properties,  can be found in \cite{Bier}.     Nevertheless,  determining whether or not an inductive limit is regular remains important, as one can see for example, in \cite[Thm. 34, p. 1371]{Dab}.   In \cite{Qiu1},  Qiu defined the concept of a quasi-locally complete space (denoted as quasi-fast complete in \cite{Qiu1}), in which each bounded set is contained in abounded set that is a Banach disk in a coarser locally convex topology, and proves that if an inductive limit of strictly webbed spaces is quasi-locally complete, then it is regular.  An important aspect of Qiu's work is that many spaces are quasi-locally complete.  Motivated by those results,  the assumption of a Banach disk is replaced with that of a Baire disk.  We generalize Qiu's theorem on the regularity of inductive limits of  strictly webbed spaces that are  quasi-locally complete \cite[Thm. 1, p. 57]{Qiu1} to that of being quasi- locally Baire. Thus, a condition for regularity is obtained that does not require any completeness on the bounded disks.    Examples are described to show that the collection of quasi-locally Baire spaces is wider than that of quasi-locally complete spaces.   \\

\noindent   Throughout this paper, {\bf space} refers to a Hausdorff  locally convex space  $E = (E, \tau)$ over $\mathbb{K} = \mathbb{R} \mbox{ or } \mathbb{C}$.  We will call a convex, balanced set a  {\bf disk}.    Given a disk $B$ of a space $E$, the linear span of $B$ is denoted by $E_{B}$.  Equipped with with the linear topology given by the Minkowski functional of $B$,  $E_{B}$ is a seminormed space.  If additionally the set $B$ is  bounded, then  with this topology, $(E_{B}, \mu_{B})$ is a normed space and we write  $(E_{B}, ||\cdot ||_{B})$.  When this normed space is complete (resp. Baire),  $B$ is called a {\bf Banach disk} (resp. {\bf Baire disk}).  A space $E$ is {\bf locally complete} (resp. {\bf locally Baire}), if every bounded subset is contained in a Banach (resp. Baire) disk, c.f. \cite[5.1.6 (iv), p. 152]{PCB} (resp. \cite[p.45]{TG1}.  Certainly, every locally complete space is locally Baire, and the existence of incomplete normed Baire spaces \cite{Klis} shows that the reverse implication does not hold.   Consider  spaces $(E_{n}, \tau_{n}), \, n\in\mathbb{N}$ such that $E_{n}\subset E_{n+1}$, with continuous injections $id : E_{n} \rightarrow E_{n+1}$  for each $n\in\mathbb{N}$.  Put  $E = \bigcup_{n\in\mathbb{N}} E_{n}$. Denote by $(E, \tau_{ind}) = indlim_{n}(E_{n})$, the corresponding locally convex inductive limit. The space  $(E, \tau_{ind})$ is {\bf regular} if given any bounded subset of $E$, there exists $m\in\mathbb{N}$ such that the subset is contained in and bounded in $(E_{m}, \tau_{m})$.   We will make use of results regarding absolutely convex webs on locally convex spaces.  The relevant references in use here are \cite[Appendix]{RR} and \cite{WR}.    Any other unspecified notation follows that of \cite{RR}.    The  main result of this paper is:  \\

\begin{thm}  Suppose each $E_{n}$ is strictly webbed. 
\begin{enumerate}
\item[(a).]   If $E = indlim(E_{n})$ is quasi-locally Baire, then $E$ is regular.   
\item[(b).]  If  each $E_{n}$ is quasi-locally Baire and $E$ is regular, then $E$ is locally Baire.  \\
\end{enumerate}
\end{thm}

\noindent  Part (a) of Theorem 1 generalizes \cite[Thm. 1, p. 57]{Qiu1}  and  part (b) generalizes \cite[Thm. 1, p. 46]{TG1}.    In \cite[p. 56]{Qiu1} the following is defined:

\begin{defn}  A space $E = (E, \tau)$ is {\bf quasi-locally complete (quasi-fast complete)}, \cite[p. 56]{Qiu1},  if given any bounded subset $A$  of $E$, there exists a   locally convex topology $\rho$ that is coarser than $\tau$ on $E$  such that $A$ is contained in a $\rho-$ Banach disk.  \\
\end{defn}

\noindent  In \cite[Thm. 1, p. 57]{Qiu1} it is proven that if each $E_{n}$ is strictly webbed  and if $E = indlim(E_{n})$ is quasi-locally complete, then $E$ is regular. We generalize quasi-local completeness and obtain a generalization of this result.  \\

\begin{defn}  A space $E = (E, \tau)$ is {\bf quasi-locally Baire} if given any bounded subset $A$  of $E$, there exists a  locally convex topology $\rho$ on $E$  coarser than $\tau$  such that $A$ is contained in a $\rho-$ Baire disk.  \\
\end{defn}

\section{Proof of Theorem 1.}  

\noindent  Let us begin with a  version of the localization theorem for strictly webbed spaces, generalizing \cite[Lemma 1, p. 56]{Qiu1}.  

\begin{lem}  Suppose $E = indlim_{n} E_{n}$ such that each $E_{n}$ is strictly webbed, and that $F$ is a metrizable Baire space.  Let  $T:F \rightarrow E$ be linear with a closed graph, and assume $A$ is any bounded subset of $F$.  Then there exists $m\in\mathbb{N}$ such that $T(A)$ is contained in and bounded in  $E_{m}$.
\end{lem}

\noindent {\it Proof}:  For each $n\in\mathbb{N}$, let $\mathcal{W}^{(n)}$ denote a strict web on $E_{n}$.  A strict web $\mathcal{W}$ on $E$ is created by defining the $k$th layer to be the collection of $k$th layers of $\mathcal{W}^{(n)}$, as $n$ varies in $\mathbb{N}$.  Thus, a collection $\{W_{k} : k\in\mathbb{N}\} \subset \mathcal{W}$ is a strand  if, and only if, it is a strand of some $\mathcal{W}^{(n)}$, \cite[p.163]{RR}.  Now let $A\subset F$ be any bounded set, and apply the localization theorem to $T$, \cite[Cor. 1, p. 722]{WR} :  There exists a strand $(W_{k})$ of $\mathcal{W}$ and scalars $\alpha_{k}$, $k\in\mathbb{N}$, such that $T(A)\subset \alpha_{k} W_{k}$.  The construction of $\mathcal{W}$ is such that, for some $m\in\mathbb{N}$, $(W_{k})$  is a strand of $\mathcal{W}^{(m)}$ in $E_{m}$.  Because strict webs are compatible with the relevant linear topologies, \cite[p. 156]{RR}, given any absolutely convex zero neighborhood $U$ in $E_{m}$, there exists $K\in\mathbb{N}$ such that $W_{K}\subset U$.  Hence, $T(A)\subset \alpha_{K} W_{K} \subset U$, proving that $T(A)$ is bounded in $E_{m}$.    \/ \/ $\Box$   \\
\bigskip

\noindent  We now prove the main result.   \\

\noindent {\it Proof of Theorem 1}:  (a):  Denote the inductive limit topology on $E$ by $\tau_{ind}$.   Let $A\subset E$ be any bounded set.  There exists a locally convex topology $\rho$, coarser than $\tau_{ind}$ and a  $\rho$ - Baire disk $B\supset A$.  That is, $(E_{B}, ||\cdot ||_{B})$ is a Baire space.    Moreover, as $B$ is $\rho$ - bounded, the normed topology of $(E_{B}, ||\cdot ||_{B})$ on the linear space $E_{B}$  is finer than the  topology inherited from $\rho$; in other words, the identity map $id: (E_{B}, ||\cdot ||_{B}) \rightarrow (E, \rho)$ is continuous, hence, has a closed graph.   Observe that $id:  (E, \rho) \rightarrow (E_{B}, \tau_{ind})$ has a closed graph.  We have before us a linear map having a closed graph from a Baire space to a strictly webbed space.  Apply the closed graph theorem \cite[Thm. 14, p. 716]{WR} to conclude  that  $id:  (E_{B}, ||\cdot ||_{B}) \rightarrow (E_{B}, \tau_{ind})$  is continuous and thus has a closed graph.  The proof is finished once we apply Lemma 1:  There exists $m\in\mathbb{N}$ such that $id(B) = B\subset E_{m}$ and $B$ is bounded in the topology of $E_{m}$.  In particular, $A\subset B$ is  bounded in $E_{m}$. \\

\noindent  (b): Assume each $E_{n}$ is quasi-locally Baire and $E$ is regular. Note that because $E$ is regular, the inductive topology is Hausdorff \cite[8.5.13 (iii), p. 286]{PCB}.   Let $A\subset E$ be any bounded set.  There exists $m\in\mathbb{N}$ such that $A\subset E_{m}$ and is bounded in the topology $\tau_{m}$ of $E_{m}$.  Next, there exists a locally convex topology $\rho_{m}$ on $E_{m}$, coarser than $\tau_{m}$, such that $A\subset B$, and $B$ is a $\rho_{m}$ - Baire disk.  Let $F$ denote the span of $B$ with the norm generated by the Minkowski functional of $B$.  The linear map $id : F \rightarrow E$ from a Baire space to a strictly webbed space has a closed graph.  Once again, we apply the closed graph theorem  \cite[Thm. 14  p. 716]{WR}, we conclude that this identity map is continuous.  Therefore, $A$ is contained in $B$ and $B$ is a Baire disk.  We conclude that $(E, \tau_{ind})$ is locally Baire.     \/ \/ $\Box$   \\

\bigskip
\section{Examples.}

\noindent  Clearly, every locally Baire space is quasi-locally Baire, and every quasi-locally complete space is quasi-locally Baire.  Below are examples that distinguish these  collections of spaces.   

\begin{ex}  Quasi-locally Baire spaces that are  not locally Baire.
\end{ex}

\noindent  In fact, such spaces exist in abundance.   By \cite[4.6.7 (iv), p. 131]{PCB}, on every  infinite dimensional Banach space $(E, ||\cdot ||)$ there exists a strictly finer norm $||\cdot||_{r}$ such that $(E, ||\cdot ||_{r})$ is not barrelled.  The space $(E, ||\cdot ||_{r})$ is quasi-locally Baire, in fact it is even quasi-locally complete.  On the other hand, as a non-barrelled space, $(E, ||\cdot ||_{r})$ certainly cannot be a Baire space and in particular, the closed unit ball of $(E, ||\cdot ||_{r})$ cannot be a Baire disk.    \/ \/ \/ $\Box$   \\

\begin{ex}  A quasi-locally Baire space that is not quasi-locally complete with respect to compatible topologies.   
\end{ex}

\noindent  By \cite{Klis}, there exists a normed space $(E, ||\cdot ||)$  that is Baire and incomplete.  This space  is locally Baire, so also quasi- locally Baire.   No locally convex topology on $E$ that is compatible with the duality of the pair $(E, E^{\prime})$ can be locally complete because Mackey's theorem \cite[Prop. 3.4.3, p. 198, Thm 3.5.3, p. 209 ]{Hor} would imply that  $(E, ||\cdot ||)$ must also be locally complete.  Such completeness would create a contradiction given that a normed space is locally complete if, and only if, it is complete.  In particular, no compatible locally convex topology coarser than the normed topology can be locally complete.   \/ \/ \/ $\Box$   \\

\noindent  Without further contemplation, Theorem 1 (b) reduces to a known situation. Indeed, a  Baire space that is strictly webbed is a Fr\'{e}chet space, by \cite[Cor. 3, p. 722]{WR}.  It behooves us to find strictly webbed spaces that are quasi - locally Baire, or at least locally Baire, and not Baire.  Such spaces are also abundant:   \\

  \begin{ex}  Strictly webbed locally Baire spaces that are not Baire.   
\end{ex}
\noindent  Consider any infinite dimensional Banach space with its weak topology, $(E, \sigma(E, E^{\prime}))$.  With its normed topology such a space is strictly webbed, so as a continuous linear image of a strictly webbed space, $(E, \sigma(E, E^{\prime}))$ is strictly webbed   \cite[Thm. 29, p. 731]{WR}.   Of course, $(E, \sigma(E, E^{\prime}))$ is not even barrelled.    Meanwhile, as closed bounded disks are the same in all compatible topologies,   $(E, \sigma(E, E^{\prime}))$ is locally complete, hence locally Baire.  \/ \/ \/ $\Box$  \\

\end{document}